\documentclass[11]{article}

\usepackage[noend]{algorithmic}
\usepackage{algorithm}
\usepackage{authblk}
\usepackage{amsmath}
\usepackage{url}
\usepackage{multirow}

\usepackage{cite}
\usepackage{nicefrac}       %
\usepackage{amsfonts}
\usepackage{amssymb}
\usepackage[skip=1pt]{caption}  %

\title{Job Shop Scheduling with Integer Programming, Shifting Bottleneck, and Decision Diagrams:\\A Computational Study\thanks{B. King and R. Hildebrand were partially funded by ONR Grant N00014-20-1-2156 and AFOSR grant FA9550-21-1-0107.  Any opinions, findings, and conclusions or recommendations expressed in this material are those of the authors and do not necessarily reflect the views of the Office of Naval Research or the Air Force Office of Scientific Research.}}

\author[1]{Brannon B. King}
\author[1]{Robert Hildebrand}
\affil[1]{Department of Computer Science, Virginia Tech}
\affil[2]{Grado Department of Industrial Systems and Engineering}

\begin{document}

\maketitle              %

\begin{abstract}
We study heuristic algorithms for job shop scheduling problems. 
 We compare classical approaches, such as the shifting bottleneck heuristic with novel strategies using decision diagrams.  Balas' local refinement is used to improve feasible solutions. Heuristic approaches are combined with Mixed Integer Programming  and Constraint Programming approaches.  We discuss our results via computational experiments.
\end{abstract}

\section{Introduction}

The job shop scheduling problem (JSP) has long been a challenging area in operations research, historically tackled through disjunctive integer programming formulations that often yield poor linear programming relaxations. Constraint programming (CP) has emerged as a more effective approach for these problems, outperforming traditional mixed-integer programming (MIP) models. Heuristics, including scheduling and dispatching rules, have been extensively studied and applied to provide feasible solutions that can be further refined through local search algorithms. This paper explores the integration of problem-specific heuristics into modern solvers, with a focus on compound heuristic approaches utilizing Decision Diagrams (DDs) and Balas's local search methods.

Scheduling problems are known for having poor linear programming relaxations. To quote \cite{martinNewApproachComputing1996a}: ``In spite of a great deal of effort, the disjunctive integer programming formulation of the job-shop problem appears to be of little assistance in solving instances of even moderate size; furthermore, its natural linear programming relaxation has been shown to give very poor lower bounds for the problem.'' Solving the problem for a subset of machines/jobs seems to be the go-to mechanism for finding a better lower bound, and the solvers rely mostly on branching to improve the bounds.

On the other hand, primal heuristics for the problem abound \cite{guptaReviewSchedulingRules1989, panwalkarSurveySchedulingRules1977, allahverdiThirdComprehensiveSurvey2015}. They are often referred to as scheduling or dispatching ``rules''. And, as these provide feasible solutions, it's recommended to follow the heuristic with a local minimization algorithm. Also, approximation algorithms exist \cite{grigorevaWorstCaseAnalysisApproximation2021}. Balas~\cite{balasMachineSequencingDisjunctive1969} gave an inspiring local search algorithm many years ago based upon the critical path in a disjunctive graph representation of a scheduling problem. This heuristic is easily conjoined to modern MIP solvers.

Recent studies comparing constraint programming (CP) to mixed-integer programming (MIP) models show that CP clearly outperforms MIP in the realm of scheduling problems. See \cite{KU2016165}, \cite{DACOL2022100249}, and \cite{mipVsCP}. The formulas in those papers are fairly standard and build on long-existing formulations \cite{panStudyIntegerProgramming1997}. We find, though, that tools/mechanisms for mixing heuristics with CP tools are lacking.
 
 Recently, \emph{Decision Diagrams} have shown to be a useful approach to some types of optimization problems. These perform a decomposition of the problem based on a sequential decision-making process. Bergman et al's book \cite{bergmanDecisionDiagramsOptimization2016} on Decision Diagrams (DDs) is the starting point for our work. We assume the reader's familiarity with said book's content. We also assume familiarity with scheduling problems and their triplet notation (as in \cite{michaell.pinedoSchedulingTheoryAlgorithms2022}). At first sight, DDs appear to be nothing more than a formal method for enumerating all possible solutions to a problem, with the detection of duplicated intermediates. However, their power is found in two separate mechanisms: 1) the restricted form that uses fixed memory to generate multiple feasible solutions, and in 2) how they can provide a relaxed form of the problem. This is done by intelligently merging nodes as the tree of solutions grows too wide. Solutions that go through one or more of these relaxed nodes provide a useful dual bound.

\subsection{Contributions}
We develop heuristic approaches based on Balas's work and on Decision Diagrams.  
Our novel Decision Diagram models for the JSP encourage minimizing stored symmetry, and thus reducing computational effort.

We evaluate the effectiveness of Restricted Decision Diagrams compared to traditional heuristics for JSPs.  Lastly, we investigate the impact of warm-starting modern solvers with a heuristic solution. We then discuss the conclusions of our computational experiments. We use GUROBI\cite{gurobi} and CPLEX\cite{cplex} for MIP solvers and also CPLEX's constraint programming solver when testing CP versions of the JSP.

\section{Formulation Background}

\subsection{Job shop scheduling}
The scheduling problem denoted as $Jm || C_{max}$ refers to a specific class of job shop scheduling problems (JSPs) characterized by the goal of minimizing the makespan across multiple machines. Formulas for it are common in existing literature, e.g \cite{panStudyIntegerProgramming1997}. Formal description follows:

Let there be a set of $n$ jobs $\{J_1, J_2, \ldots, J_n\}$ and a set of $m$ machines $\{M_1, M_2, \ldots, M_m\}$. Each job $J_i$ consists of a sequence of operations $\{O_{i1}, O_{i2}, \ldots, O_{im}\}$, where each operation $O_{ij}$ must be processed on a specific machine $M_{\pi_i(j)}$ for some permutation map $\pi_i\colon [m]\to [m]$ and for a predetermined duration $p_{ij}$ without interruption. Each machine can process only one operation at a time, and each operation can be processed on exactly one machine as per its predefined sequence in a job. The objective is to find a schedule, i.e., an allocation of operations to time slots on each machine, that minimizes the makespan $C_{max}$, which is the time when the last job completes processing.

\subsubsection{Mathematical Optimization Formulation}

We will focus on the disjunctive MIP formulation, which generally solves the quickest using MIP solvers \cite{panStudyIntegerProgramming1997}. The idea is to model the problem using binary variables to represent the sequencing decisions between operations on the same machine. Here's a step-by-step formulation:

\textbf{Variables:}

\begin{itemize}
    \item $S_{ij}$: Nonnegative start time of operation $O_{ij}$.
    \item $C_{ij}$: Completion time of operation $O_{ij}$ \\ and easily collapsed into $S_{ij} + p_{ij}$.
    \item $C_{\max}$: Maximum completion time (makespan).
    \item $x_{ijkl}$: Binary variable indicating operation $O_{ij}$ follows $O_{kl}$, only existing if both are on the same machine.
\end{itemize}

\textbf{Constraints:}

\begin{itemize}
    \item Precedence Constraints: Ensure the correct order of operations within each job.
    \item Disjunctive Constraints (Eqs. 1c, 1d): Ensure that no two operations on the same machine overlap by enforcing that one must precede the other.
\end{itemize}

\textbf{Mixed Integer Programming Model:}
\begin{subequations}
\begin{align}
\min & ~C_{\max} \quad \text{subject to} \\
& S_{ij} \geq C_{i(j-1)} \quad \forall i \in [n], \forall j \in J_i, j > 1 \\
& S_{ij} \geq C_{kl} - \overline{M}(1 - x_{ijkl}) \quad \forall ij, kl : \pi_i(j) = \pi_k(l)\\
& S_{kl} \geq C_{ij} - \overline{M}x_{ijkl} \quad \forall ij, kl : \pi_i(j) = \pi_k(l)\\
& S_{ij} \geq 0 \quad \forall i, j \\
& C_{ij} = S_{ij} + p_{ij} \quad \forall i, j \\
& C_{\max} \geq C_{ij} \quad \forall i, j : O_{ij} \text{ is final task of job } i \\
& x_{ijkl} \in \{0, 1\} \quad \forall i, j, k, l : \pi_i(j) = \pi_k(l)
\end{align}
\end{subequations}
where $\overline{M}$ is a big-M multiplier, typically set to the one plus the sum of the delays: $1+\sum_{ij}p_{ij}$. Those big-M constraints may also be formulated using the solver's indicator constraint functionality.

\subsubsection{Constraint programming formulation}

See the full explanation in \cite{mipVsCP}.

\textbf{Variables:}

\begin{itemize}
    \item $I_{ij}$: Interval variable containing the start and end of operation $O_{ij}$ with width as the  duration$=p_{ij}$.
\end{itemize}

\textbf{CP Formula:}
\begin{subequations}
\begin{align}
\min~ & \max_{i\in [n]}(\text{end\_of}((I_{im})))
 \quad \text{subject to:} \\
& \text{no\_overlap}(\{I_{ij} \mid i \in [n] \}) \quad \forall j \in [m] \\
& \text{end\_before\_start}(I_{ij}, I_{i(j+1)}) \quad \forall i \in [n], j \in [m-1]\\
& I_{ij} = \mathrm{IntervalVar}(p_{ij}) \quad \forall i\in [n], j \in [m]
\end{align}
\end{subequations}

\subsection{Existing Work regarding scheduling via DD}

Bergman et al. \cite{bergmanDecisionDiagramsOptimization2016} discuss how a DD can solve scheduling problems in general and give this example: the single-machine makespan minimization with sequence-dependent delays: $1|p_{ij}|C_{max}$. They don't use a binary expansion tree; instead, they represent the problem as a permutation of possible orderings, known as a multivalue decision diagram (MDD), as shown in Table~\ref{table:I}.

\begin{center}
\captionof{table}{Bergman's simple state operators}
\label{table:I}
\noindent\begin{tabular}{r|p{2.575in}}
    \hline
    State & $S_j$ holds the $j$ jobs already done \\
    Transition & $S_j \cup {x}$ where $x \in [n]$ and $x \notin S_j$ \\
    Cost & the delay from $S_j$ to $S_{j+1}$ via $x$ \\
    \hline
\end{tabular}
\end{center}

They do not give merge and split definitions for their scheduling example. They do give merge operations for other problems, namely, maximum independent set, maximum cut, and maximum 2-SAT. They also discuss how some merge operations might be possible for a scheduling problem if we separate the jobs already done into two groups: those that are surely done no matter what path arrives at a given state, and those that are done in at least one path arriving at a given state. This latter group is the "maybes". In \cite{bergmanDecisionDiagramsOptimization2016}'s terms, the two groups are ``All'' and ``Some''. The maybes doesn't exist in a full expansion because, in that context, we don't have any nodes that represent more than one unique solution.

Hooker, in~\cite{hookerJobSequencingBounds2017}, expands on the above ideas for a more complicated scheduling problem shown in Table~\ref{table:II}, minimizing tardiness given release dates and due dates $1|r_j,d_j|\sum T_j$.

\begin{center}
\captionof{table}{Hooker's tardiness operations}
\label{table:II}
\noindent\begin{tabular}{r|p{2.575in}}
    \hline
    State & a tuple $S_j=(V, U, f)$:\newline
    \hspace*{2pt} $V$ holds up to $j$ jobs surely done,\newline
    \hspace*{2pt} U holds jobs done on some route,\newline
    \hspace*{2pt} $f$ is the running completion time \\
    Transition & ($V \cup x, U \cap x, \max(r_x, f) + p_x)$ \newline
    \hspace*{2pt} where $x \in [n]$ and $x \notin V$ \\
    Cost & the delay from $S_j$ via $x$, the value of $p_x$ \\
    Merge & ($V \cap V', U \cup U' , \min(f, f'))$\\
    \hline
\end{tabular}
\end{center}

His paper demonstrates that these operations are sufficient to ensure that the relaxed tree contains a path that represents a dual (in this case, lower) bound. He also gives a mechanism for proving any merge operation to be sufficient. He expands his effort with a later paper, \cite{hookerImprovedJobSequencing2019}, where he includes merge operations for sequencing with time windows, time-dependent processing times, sequence-dependent processing times, and state-dependent processing times. They all follow a pattern very similar to the one given above.

Moreover, in \cite{hookerImprovedJobSequencing2019}, Hooker suggests optimizing the Lagrangian dual where he incorporates an additional penalty for sequences that repeat operations -- a common infeasibility in a relaxed DD representation of a schedule. It was inspired by \cite{bergmanLagrangianBoundsDecision2015a}. The paper also includes a table of timed runs on the CPW and Biskup-Feldman public datasets containing tardiness problems. (We were unable to make this succeed in our context. Perhaps it was due to some failure of our model to meet the necessary assumptions, or it required an extremely high number of iterations for convergence, or coding error.)

Building on Hooker's work, \cite{vandenbogaerdtMultimachineSchedulingLower2018} provides two recent papers tackling multi-machine scheduling. In \cite{vandenbogaerdtMultimachineSchedulingLower2018}, they focus on a tardiness problem with substantial state; $(V, U, f, t, f^u, t^u, g)$. $f, t, f^u, t^u$ are all vectors, where $f$ refers to running completion times per machine, $t$ refers to running release times, and the superscript $u$ implies the same for the maybes (the items included in some ancestral lines but not all). They show that their merge operation is suitable using Hooker's rules. It's given here:
\begin{align*}
(V \cap V',~U \cup U',~\min(f, f'),~\min(t, t'),\\
\quad \max(f^u, f^{u\prime}),~\max(t^u, t^{u\prime}),~\min(g, g'))
\end{align*}

In their later work, \cite{vandenbogaerdtLowerBoundsUniform2019}, they build rules for uniform scheduling over total tardiness, or $Um|r_j,st_{jj'},d_j|\sum T_j$. They track the current machine as part of the state. This leads to a notable limitation; their merge operation is only allowed to merge nodes where the current machine matches. They generally follow the patterns given above for tardiness problems.

\subsection{Decision Diagrams}
See~\cite{castroDecisionDiagramsDiscrete2022} for a recent review of decision diagrams (DDs) for discrete optimization.  We follow that survey for some details.

\section{DD Operators for the JSP}
In this section we present several different models that store and transition state in DDs. All solve the JSP, but not at the same efficiency. They are not the only possible models. For merging state, additional information must be stored, and that is covered in a separate section.

If we assign each operation a unique identifier, we have all feasible solutions as the permutations of those identifiers. Of course, most permutations are infeasible in that some job's operations may be out of order. Moreover, if we have any permutation, either complete or partial, we can trivially compute its completion times, feasibility, and, for partial sets, a lower bound on the completion max. See Algorithm \ref{alg:cfp1}.

\begin{algorithm}[H]
\caption{Cost from partial solution (CFP)}\label{alg:cfp1}
\begin{algorithmic}[2]
\REQUIRE $X$ is a list of to-be-done operations.
\REQUIRE $\textbf{f}^M$ is machine completion times, \textbf{0} by default.
\REQUIRE $O$ operations already done with times $\textbf{f}^O$.
\FOR{$x \in X$}
    \STATE $s \gets \textbf{f}^M_{mach(x)}$
    \IF{$\exists ~pre(x)$}
        \IF{$pre(x) \notin O$}
            \RETURN Error: Missing Prerequisite
        \ENDIF
        \STATE $s \gets \textsc{max}(s, \textbf{f}^O_{pre(x)})$
    \ENDIF
    \STATE $\textbf{f}^O_x \gets s + \textsc{delay}(x)$
    \STATE $\textbf{f}^M_{mach(x)} \gets s + \textsc{delay}(x)$
\ENDFOR
\RETURN{$\textsc{max}(\textbf{f}^M)$}
\end{algorithmic}
\end{algorithm}

Given a partial ordering and the next operation to be concatenated to it, we can trivially compute the change in completion times brought about by the additional operation. We can also update the running $C_{max}$ if it is changed by this additional operation. 

\subsection{The basic permutation model}

Since we're just storing a partial ordering, all we need is a list. We make use of some helpers such as $pre(x)$, which returns the operation required right before $x$ if there is one. $mach(x)$ returns the required machine for operation $x$. $delay(x)$ is the delay for operation $x$, commonly referred to as $p_{ij}$. $trailer(x)$ contains the sum of the operation times that must follow operation $x$, where this could be either the amount of work left on the job of $x$ or the machine required for $x$ (or the maximum of both of those). Given a list of operations with labels in $[n]$, we get Table~\ref{table:III}.

\begin{center}
\captionof{table}{JSP-for-MDD Model 0}
\label{table:III}
\noindent\begin{tabular}{r|p{2.575in}}
    \hline
    State & $S :=V$, an ordered list of ops. done so far \\
    Cost & $c(S,x) := \textsc{CFP}([V, x])$ \\
    Transition & $\phi(S,x) := [V, x]$ : $x \in [n] \backslash V$, $pre(x) \in V$ \\
    \hline
\end{tabular}
\end{center}

If we make our model slightly more advanced, we can cache the completion times for each operation (in $\mathbf{f}^O$) and each machine (in $\mathbf{f}^M)$, we get Table~\ref{table:IV}.

\begin{center}
\captionof{table}{JSP-for-MDD Model 1}
\label{table:IV}
\noindent\begin{tabular}{r|p{2.575in}}
    \hline
    State & a tuple $S := (V, \mathbf{f}^O, \mathbf{f}^M)$ \\
    Cost & $c(S,x) := \max(\mathbf{f}^M_{mach(x)}, \mathbf{f}^O_{pre(x)}) + delay(x)$ \\
    Transition & $\phi(S,x) := (V \cup \{x\}, \newline ~~~~~\mathbf{f}^O_x\gets c(S,x), \mathbf{f}^M_{mach(x)}\gets c(S,x))$ \\
    \hline
\end{tabular}
\end{center}

It's implied that $(V, \mathbf{f}^O, \mathbf{f}^M)$ are independent for every state -- copied from the parent state and then modified/extended. When storing the states in each layer, it's useful to store duplicate states only once. When comparing these states, all fields of the tuple must be compared. Notice that only states within a given layer will have equal cardinality for $V$ (unless states are merged in some way that violates that, as discussed in the merge section below).

\subsection{A model to detect more symmetry}

When growing a tree of possible solutions, such as is done by decision diagrams, one may arrive at equivalent solutions through differing paths. In this context such solutions are symmetrical. With a goal of keeping the tree of possible solutions as small as possible, DDs benefit from any reduction in state space. We want to avoid symmetry in the storage of our accumulated state.

The above model is not bad, but notice (Table~\ref{table:V}) that many of the stored completion times can never eclipse any machine's finish time. Sometimes those may happen in some other order and still produce the same state; we redesign the state to capture that symmetry.

\begin{center}
\captionof{table}{JSP-for-MDD Model 2}
\label{table:V}
\noindent\begin{tabular}{r|p{2.575in}}
    \hline
    State & a tuple $S := (V, V_L, \mathbf{f}^O, \mathbf{f}^M)$; \\
    Cost & $c(S,x) := \max(\mathbf{f}^M_{mach(x)}, \mathbf{f}^O_{pre(x)}) + delay(x)$ \\
    Transition & $\phi(S,x) := (V \cup \{x\} \backslash pre(x), V_L \cup pre(x), \newline ~~~~~\mathbf{f}^O_x\gets c, \mathbf{f}^M_{mach(x)}\gets c)$\\ 
    \hline
\end{tabular}
\end{center}

 $V_L$ refers to those operations that are long-done, those that will never need to be used as an immediate prerequisite. \emph{For hash and comparison we ignore the completion times of operations in $V_L$.} If $pre(x) \notin \textbf{f}^O$ you can return zero, although tracking some completion time for items in $V_L$ can be handy for generating the final schedule at the end. 
 
 In Table~\ref{table:VI}, we show the number of nodes for full expansion averaged over 10 random problems to demonstrate that the number of nodes expanded is reduced by capturing more symmetry:

\begin{center}
\captionof{table}{Node expansion demonstration}
\label{table:VI}
\noindent\begin{tabular}{r|l|l}
    & \textbf{Model 1} & \textbf{Model 2}\\
    \hline
    Valid nodes, 4x5 JSP & 1245k & 793k \\
    Duplicates & 1017k & 776k \\
    \hline
    Valid nodes, 3x10 JSP & 726k & 528k \\
    Duplicates & 638k & 504k \\
    \hline
\end{tabular}
\end{center}

\subsection{Modeling disjunctives directly}

While the MILP formulation described above relies on real variables for start and completion times, the values for these variables can be fully determined from a fixed set of binary variables, denoted as \( x \). This is common in scheduling problems. Consequently, this can be directly modeled as a binary decision diagram (BDD), as opposed to the multivalue or MDDs mentioned earlier. For each disjunction, a binary variable \( x_i \) determines whether the path is forward or reverse. In this context, forward implies that \( i < j \) in an adjacency matrix representation of the disjunctive graph of the problem, as explained in \cite{balasMachineSequencingDisjunctive1969}. Following the approach in that paper, it is assumed that the process begins with a feasible set of disjunctives, all oriented forward.

This approach encompasses a significantly larger number of variables and consequently necessitates many more layers compared to MDDs. For example, the well-known JSP problem \emph{abz5} involves 100 operations but 900 disjunctives. (The number of dijunctives matches the number of binary variables shown in the Gurobi log, but it can also be easily computed from the number of bidirectional arcs in 10 cliques of 10 nodes each: $10n(n-1)$.) Grouping these bits into bytes could reduce the number of layers by a factor of eight while simultaneously increasing the width of each layer by the same factor. Additionally, this method might potentially reduce the 256 possible expansions by eliminating common 3- or 4-cycle problems from the possible values, although this approach has not yet been explored.

Given a disjunctive model, including invalid models with cycles, we can always find the longest path through that model in polynomial time or less. That's done trivially with a flow LP or an adjacency matrix selection LP. The challenge is to reduce the number of LP calls needed, as one for each state (or every other state as shown below) gets expensive. This approach was generally unhelpful although it did work correctly; it's included here in Table~\ref{table:VII} for contrast.

\begin{center}
\captionof{table}{Basic BDD-with-LP for layer $j$}
\label{table:VII}
\noindent\begin{tabular}{r|p{2.575in}}
    \hline
    State & $V$ holds the reversed disjunctives' index \\
    Transition & $V'=[V, j]$ if $x=1$ else $V'=V$ \\
    Cost & LP(V) \\
    \hline
\end{tabular}
\end{center}

With this approach, each node is unique; there is no duplication and no reduction in possible values as you progress through the layers. You cannot cull states that have cycles in their graph as later reversals may eliminate those cycles. As with all DD approaches, it would be worth it to cull via some known primal bound, heuristically determined.

\subsection{Merging state}

Both Model 1 and Model 2 support merge operations. Following \cite{hookerJobSequencingBounds2017} and \cite{vandenbogaerdtMultimachineSchedulingLower2018}, we add additional fields to the state: $V_s$ and $\textbf{f}^s$, where $s$ stands for ``some'', meaning that some ancestral path covered the items in this set. For merging Model 1 state, $S' \gets S \oplus \overline{S}$:
\begin{subequations}
\begin{align}
    \textbf{f}^{M\prime}_i & \gets \min(\textbf{f}^M_i, \overline{\textbf{f}^M_i}) \quad \forall i \in M\\
    V' & \gets V \cap \overline{V}\\
    \textbf{f}^{O\prime}_x & \gets \max(\textbf{f}^O_x, \overline{\textbf{f}^O_x}) \quad \forall x \in V'\\
    V_s' & \gets (V \cup \overline{V}) - V'\\ 
    \textbf{f}^{s\prime}_x & \gets \max(\textbf{f}^O_x, \overline{\textbf{f}^O_x}, \textbf{f}^s_x, \overline{\textbf{f}^s_x}) \quad \forall x \in V_s'
\end{align}
\end{subequations}

Hooker \cite{hookerJobSequencingBounds2017} gives two criteria for a valid merge: any possible control values leaving $S'$ must be a valid relaxation of the same control leaving $S$, and both $S$ and $\overline{S}$ must be relaxed/interchangable. The latter criteria holds from the lack of order-dependent operations in the above formulation. The former criteria holds in that we always take the minimum score for our $C_{max}=\max(\textbf{f}^{M\prime}_i)$, and we always use the worst-case prerequisite completion time when pulling items from $\textbf{f}^O$. Model 2 is merged similarly with $V_L$ receiving equivalent treatment to $V$.

We used the merge operation to validate that Bergman's branch-and-bound algorithm \cite{bergmanDecisionDiagramsOptimization2016} worked for JSP. It did not scale well, though. We were unable to make the relaxed DD give a better lower bound than the linear polytope of the disjunctive formulation, so we do not include this in our computational results below. As part of that, we pruned nodes that exceeded the current primal bound, which is known as LocB pruning in this context \cite{gillardLocB}. See the other thoughts on the algorithm in the appendix.

\section{A simple local search refinement}

In Balas' original paper \cite{balasMachineSequencingDisjunctive1969} about solving JSP via disjunctive graph iteration, he included this special proposition:

\textbf{Proposition of Balas 1969:} 
\begin{quote}
Let $C_h$ be a critical path in $G_h$ [which is a DAG]. Any graph $G_k$ obtained from $G_h$ by complementing one arc $(i, j) \in C_h$ is circuit-free.
\end{quote}

\textbf{Balas' proof:} We know that $(i,j)$ must be the longest path from $i$ to $j$ as it is part of $C_h$. However, it is also the shortest path from $i$ to $j$ or we would have chosen the longer path. Because $i$ to $j$ is the only path, we can reverse it without creating a circuit.

That leads us to this refinement algorithm as a local search method (Alg. \ref{alg:lns1}):
\begin{algorithm}[H]
\caption{Balas Local Refinement 1 (LNS1)}\label{alg:lns1}
\begin{algorithmic}[2]
\REQUIRE $g = (V, \mathcal A)$ is a DAG of the fixed conjunctive arcs. 
\REQUIRE $W_e \geq 0$ for $e \in \mathcal A$ as the weight of each arc.
\REQUIRE $s_{parent}$ is the $s$ value from the caller.
\STATE $p \gets \textsc{longest\_path}(g)$
\STATE $s \gets \sum_{(u,v) \in p} W_{uv}$
\IF{$s > s_{parent}$}
    \RETURN $s$ \hfill //~remove this to search more space
\ENDIF
\FOR{$e \in p$}
    \IF{not $e.fixed$}
        \STATE $g \gets \textsc{remove\_arc}(g, e)$
        \STATE $g \gets \textsc{add\_fixed\_arc}(g, e.v, e.u, W_{vu})$
        \STATE $t \gets \textsc{recurse}(g, W, s)$
        \IF{$t < s$}
            \STATE $s=t$ \hfill //~can also track swaps here
        \ENDIF
        \STATE $g \gets \textsc{remove\_arc}(g, e.v, e.u)$
        \STATE $g \gets \textsc{add\_fixed\_arc}(g, e)$
    \ENDIF
\ENDFOR
\RETURN{$s$} \hfill //~return swaps also if desired
\end{algorithmic}
\end{algorithm}

\section{Computational Experiments}

\subsection{Use as a Heuristic}

In this section we show how a restricted Model 2 compares to other common heuristics including MOR, MWR, and the shifting bottleneck \cite{adamsShiftingBottleneckProcedure1988a}. We include multiple widths for the restricted DD, but this parameter does not substantially improve the bound it computes; see \cite{hookerJobSequencingBounds2017}.

\textbf{Observation:} The restricted DD always produces at least one feasible solution. Given any partial solution that is feasible, the remaining operations can always be added in a feasible order. Note that this does not hold if you filter the nodes in the DD with anything additional to the maximum width filter. For example, if you filter nodes whose $C_{max}$ exceeds some threshold (in addition to reducing row width), you may filter all possible paths toward the conclusion of the DD.

In Table~\ref{table:VIII}, we run each heuristic over twenty random $10\times 10$ JSP instances, making use of Gurobi \cite{gurobi}, CPLEX \cite{cplex}, and Job Shop Library \cite{jobShopLib}. $\mathcal{W}$ refers to the maximum width of the DD. Overage is how far the final bound was above the optimum. We use default solver settings with the exception of Gurobi's AggFill=10 and GomoryPasses=1, which was recommended by their tuner for these problems. The random instances are integer and similar to those published with \cite{adamsShiftingBottleneckProcedure1988a}).

\begin{center}
\captionof{table}{Avg. heuristic overages for 20 random}
\label{table:VIII}
\begin{tabular}{r|l|l|l}
    \hline
    \textbf{Heuristic} & \textbf{Time} & \textbf{Overage} & \textbf{After LNS1} \\
    \hline
    Gurobi, MIPGap=0.25 & 0.27s & 11.1\% & 9.18\% \\
    Shifting Bottleneck & 4.0s & 15.8\% & 15.3\% \\
    Restricted DD, $\mathcal{W}$=200 & 0.25s & 18.5\% & 14.6\% \\
    Restricted DD, $\mathcal{W}$=400 & 0.50s & 16.8\% & 12.5\% \\
    \hline
    \multicolumn{2}{r|}{Most Work Remaining (MWR)} & 26.9\% & 16.1\% \\
    \multicolumn{2}{r|}{Most Ops. Remaining (MOR)} & 20.1\% & 14.1\% \\
    \multicolumn{2}{r|}{Shortest Proc. Time (SPT)} & 90.0\% & 40.1\% \\
    \hline
\end{tabular}
\end{center}

\begin{center}
\captionof{table}{Avg. heuristic overages for 18 from JSPLIB}
\label{table:IX}
\noindent\begin{tabular}{r|l|l|l}
    \hline
    \textbf{Heuristic} &\textbf{Time} & \textbf{Overage} & \textbf{After LNS1} \\
    \hline
    Gurobi, MIPGap=0.25 & 1.6s & 7.77\% & 5.85\% \\
    Gurobi, MIPGap=0.40 & 0.12s & 25.8\% & 12.5\% \\
    \hline
    \multicolumn{2}{r|}{Shifting Bottleneck} & 22.8\% & 21.9\% \\
    \multicolumn{2}{r|}{Restricted DD, $\mathcal{W}$=200} & 14.4\% & 11.9\% \\
    \multicolumn{2}{r|}{Restricted DD, $\mathcal{W}$=400} & 11.6\% & 9.70\% \\
    \multicolumn{2}{r|}{Most Work Remaining (MWR)} & 31.4\% & 19.8\% \\
    \multicolumn{2}{r|}{Most Ops. Remaining (MOR)} & 29.9\% & 23.6\% \\
    \multicolumn{2}{r|}{Shortest Proc. Time (SPT)} & 80.8\% & 39.7\% \\
    \hline
\end{tabular}
\end{center}

Interestingly, (Table~\ref{table:IX}) the DD approach worked better on the real-world problems found in JSPLIB \cite{jsplib} -- it's eighteen $10\times 10$ problems. Similarly, the small 0.25 MIPGAP for Gurobi was much more difficult to achieve on the JSPLIB problems.

The comparison is a little bit unfair, in that the Restricted DDs generate many feasible solutions for the LNS1 whereas the top four only produce a single solution to be refined. However, it shows that the local refinement eliminates the need for the shifting bottleneck heuristic. 

We recognize that there are more sophisticated versions of the shifting bottleneck algorithm, e.g. \cite{balasJobShopScheduling2008a}. There are also a variety of other local search mechanisms designed for JSP that are far more sophisticated and far-reaching, typically built on taboo search, e.g. \cite{taillardParallelTabooSearch1994}, \cite{nowickiAdvancedTabuSearch2005}. We did not consider simulated annealing nor any evolutionary algorithm as part of this research either, though papers on those approaches for JSP abound.

\textbf{Regarding Gurobi NoRel:} We ran Gurobi's NoRel heuristic for 4 seconds on the same 20 problems. It failed to find any solutions on 80\% of the problems, but on the other four, it found solutions within 3\% of optimal. Note that Gurobi can solve a 10x10 disjunctive program in 2 to 8 seconds on our test machine, so running a 4-second heuristic for it would not make sense generally. The NoRel runtime has to be specified as an input.

\textbf{Regarding runtimes:} Note that the Shifting Bottleneck (SB) without readjustment of machines in $M_0$ plus the LNS at the end still achieves 18\%. This takes about two seconds to run whereas the other takes 4 seconds per 10x10. SB can be modified to solve subproblems in parallel, which was not a part of our implementation. We rely on Gurobi to solve the $1|r_j|L_{max}$ subproblems in the SB. This takes up 90\% of the runtime for it. Generally, though, when using SB one would use Carlier's approximation \cite{carlier1982one} for the $L_{max}$ instead of solving it via a MIP solver (and probably still run them in parallel). See other ideas here: \cite{grigorevaWorstCaseAnalysisApproximation2021}, \cite{sinaiMinimizingTardinessScheduling2021}. With our DD written in Go, the 10x10 on a max-width of 400 takes half a second to run and half that time for the 200 width (using no parallelism). Most of that time is in comparison to previous nodes on the same row. Better hashing would improve that time. The LNS1 adds some additional time to that as it runs on each resulting node. This is not included in the time measurement listed above. The dispatching rule approximations obviously use a trivial amount of time.

Most of the items arriving at the bottom row of a DD tend to be similar, which comes from the sort-and-truncate approach. It needs some other selection criteria toward the top of the tree so that it keeps more diversity early on, which should give it better chances of enabling a good/unique neighborhood. The DD is quite sensitive to changes in the running $C_{max}$ computation. For example, you can use the $C_{max}$-so-far or you can add to the trailer for the remaining items to be done on each machine or you can add the work remaining on the job. Those selections all change results quite drastically. The run recorded above does not add a trailer, as going without seemed slightly better on average.

\subsection{The value of a starting point}

Here we demonstrate the value of using a heuristic to select a starting solution when solving the JSP to its optimum. In Table~\ref{table:X} we run each solver over the same 20 instances used above but give it a single starting point. The starting points are derived using MOR followed by the refinement of the LNS1 algorithm described above.

\begin{center}
\captionof{table}{Solve time with single warmstart}
\label{table:X}
\noindent\begin{tabular}{r|r|r}
    \hline
    & & \\
    \textbf{Solver on 10x10} & \textbf{Time/problem} & \textbf{With warmstart} \\
    \hline
    Gurobi, big-M & 2.1s & 1.9s  \\
    Gurobi, indicator & 4.0s & 3.6s \\
    CPLEX MP, big-M & 3.4s & 3.8s \\
    CPLEX MP, indicator & 130s & 110s \\
    CPLEX CP & 1.17s & 1.16s \\
    \hline
    & & \\
    \textbf{Solver on 12x12} & \textbf{Time/problem} & \textbf{With warmstart} \\
    \hline
    Gurobi, big-M & 36s & 43s  \\
    CPLEX MP, big-M & 66s & 60s \\
    CPLEX CP & 5.3s & 5.5s \\
    \hline
\end{tabular}
\end{center}

We conclude from this table that you should be using a constraint solver for exact solutions on this, and that the big-M path is more optimized than the indicator constraint feature, and that warm-starting it is unnecessary. Runs were made with Gurobi 11.0.2, CPLEX 22.1.1 on a i7-8750H processor.

Gurobi supports a heuristic parameter for controlling the percentage of time spent in heuristics. The default is 5\%. We explored other values from 0\% through 50\% but could find no other value to improve the average time. Increasing the parameter by 5\% generally added 5\% to the runtime.

\subsection{The value of LNS via callback}

We demonstrate (Table~\ref{table:XI}) the value of calling LNS1 in the MIP node callback (CB). We take the ordering given by the start time variables ($S$) and run the local search on that. We can actually just submit the solution given by the ordering instead of running a local optimizer on it -- herein called ``Nearest''.

\begin{center}
\captionof{table}{LNS1 in MIP-node callback}
\label{table:XI}
\begin{tabular}{r||r||r|r||r|r}
    \hline
    \textbf{Solver} & \textbf{Default} & \textbf{LNS1} & \textbf{Wins} & \textbf{Nearest} & \textbf{Wins} \\
    \hline
    Gurobi & 2.1s & 2.6s & 5.6 & 2.1s & 6.7 \\
    \hline
\end{tabular}
\end{center}

Note that we subtract the time of the callback itself, in that it is assumed that we can come up with more efficient implementations of it or curtail the calls to it as it becomes unlikely to assist. This optimization of how often the heuristic is called is a separate issue. Writing our LNS algorithm in a heuristic form that is fast enough to justify its use is nontrivial.

CPLEX, as documented, supports heuristic callbacks. However, in attempting to use them with version 22.1.1 on Linux, accessed via the Python docplex wrapper, we were unable to determine whether or not the solver was utilizing the provided heuristic solutions. No errors were given, but the incumbent scores were not updating as expected, so we did not include the numbers for it. We also have interesting numbers for FICO Xpress, but their license prohibits publishing.

\subsection{Using our state model for A*}

Unfortunately, Model 2 alone does not seem to be sufficient to allow solving a 10x10 via A*. We note that there are other efforts to make A* utilize relaxation features of the DD approach such as \cite{HORN2021105125}. The nearness of the running $C_{max}$ to reality is of critical importance in A*-search. It is possible to improve the trailer estimate by solving 2 of ten machines: see \cite{JURISCH1995145}. This is fairly quick, especially as the problems progress and most tasks have release times available. However, empirically, it's not enough accuracy to make the A* approach feasible.

\section{Conclusion}

Conclusions from our experimentation:
\begin{enumerate}
    \item Relaxed decision diagrams are useful as a simple JSP heuristic. They are not difficult to write/use and run fast. They produce better results than other common (simple) heuristics on real-world problems.
    \item Passing a start point to the solver is not useful at 14\% away from optima. Perhaps it would be worth it if you were using some more sophisticated heuristic that could generate starting points within just a few percentage points from optima.
    \item Running Balas' critical path refining, the LNS1, does generally improve a given feasible solution. It is fast to run and generally worth it.
    \item The big-M handling in Gurobi and CPLEX is superior to their indicator constraint handling at present. It may be that our $\overline{M}$ value was small enough to tip that balance.
    \item For problems where feasible solutions are rare, it may be helpful to find a nearby feasible solution in the callback if it can be done quickly. This computation is very fast for the JSP. It did help significantly on some of the test problems. Especially consider it if you use FICO Xpress.
\end{enumerate}

Ideas for future work:
\begin{enumerate}
    \item The selection of keeper nodes in the restricted DD is of critical importance. Using a basic rule like keeping the smallest 200 is a general failure -- most of the nodes with the optimum are weeded out early on. That's the curse of these scheduling problems -- the conflicts on the attractive solutions don't show up until late in the game. If we had some kind of machine learning approach that could identify bad nodes early on, we would have a higher chance of retaining the good nodes (or vice-versa). Node selection ideas from modern solvers may also apply \cite{Achterberg2009}.
    \item Another idea for retaining nodes is to try to keep a diverse set using some kind of diversity measure that would increase the likelihood of keeping the optimum path.
    \item Relaxed decision diagrams produce many infeasible solutions. Order them and you can expand these infeasible paths until you arrive at the first and best feasible solution, a best first approach similar to A*. We attempted this. However, there are so many infeasible solutions to weed through that this is generally not a viable approach for problems at scale. If we had some equivalent to cuts-for-LP, perhaps we could cut out subtrees in a way that allowed us to arrive at the best solution much faster.
\end{enumerate}

\section*{Appendix}

\subsection{Notes on existing algorithms for exact solutions via DD}

\cite{bergmanDecisionDiagramsOptimization2016} presents two general algorithms for finding an exact solution to any program representable by a DD. The first mechanism is what they term ``compiling DDs by separation'', condensed form of the algorithm given in \cite{hadzicApproximateCompilationConstraints2008}. 

Algorithm summary: Begin with a relaxed decision diagram (DD) and identify the optimal path through it. If this path violates any constraints, separate the relaxed nodes on that path into two or more replacements. Adjust the inputs feeding into the violating node so that some go to each replacement. Similarly, replicate the outputs of the violating node to each replacement. Continue this process until the optimal path is feasible.

The process requires a split or "separation" operation, which essentially undoes the merge operation, though the necessary bookkeeping for this may be expensive. If no split operation is available, a possible solution is to backtrack to the parents of the merged nodes and regenerate their children. Additionally, we assume that the arcs store and maintain the variable value (also known as the control) and the cost of traversing them. This assumption differs from our previous experiments, where we kept only the fringe nodes with the running cost in the state. Furthermore, decision diagram (DD) creation generally employs node reduction (combining identical nodes), and this reduction must be maintained after adding additional nodes to the graph. If all identical nodes are on the same layer, the check is reasonable; otherwise, it becomes too expensive.

The second algorithm represents a branch-and-bound (BnB) approach, where you branch on a cutset of exact nodes, making a new subtree pair, both relaxed and restricted, for the decedents of each node in the cutset. Cutset refers to a set of exact nodes such that any path through the tree goes through one and only one of these nodes (before hitting any relaxed nodes).

Algorithm summary: While there are nodes in the queue, remove node \(u\). Update the primal bound, which is the cost to node \(u\) plus a heuristic from it to the end, potentially determined by a restricted tree. Construct a relaxed decision diagram (DD) with \(u\) as its root. If the best relaxation is worse than the primal bound, exclude the entire \(u\) subtree. Otherwise, add the exact cutset of \(u\) to your queue and repeat the process.

\noindent Some general notes on this algorithm:
\begin{enumerate}
    \item The processing of these subtrees is parallelizable (as noted in the reference).
    \item It does not require a split operation, although it does require a working merge operation for building the relaxed trees.
    \item It doesn't require a restricted tree if you have some other heuristic mechanism that completes partial solutions, as that may also provide a reasonable primal bound, especially if it's refined by a fast local search as the final step of the heuristic. 
    \item It doesn't make any use of the dual bound for subtree exclusion. This is its fundamental weakness.
    \item Empirically, it's highly unlikely that to be able to exclude the whole relaxed tree based on its best node being worse than the current overall primal bound. Hence, you can simply return the cutset as soon as it is discovered. This eliminates the need for a merge operation.
\end{enumerate}

The algorithm makes use of two things from the relaxed DD: its best path cost and its exact cutset. \cite{bergmanDecisionDiagramsOptimization2016} gives three algorithms for selecting the cutset: take the first layer, take the last layer before any nodes are merged, or take the ``frontier'', meaning all the exact nodes that have at least one relaxed child. From that, we make these observations:
\begin{enumerate}
    \item If we merge many nodes into one, that node has a high likelihood of being very relaxed. Thus we will keep it, as it has a good score, which will in turn lead to the best path through the relaxed tree being a poor estimate of reality. Hence, again we will keep that tree's cutset, as our primal bound won't be able to exclude it.
    \item If we take some layer before we merge any nodes, our cutset will be very shallow. Shallow nodes have lower likelihood of being excluded by constraints, assuming most constraints incorporate more than the first few variables. Moreover, it is utilizing less of our DD.
    \item If we choose merge a lot in hopes of not over-relaxing any one path, we will force our cutset to be more shallow, thus getting less advantage from our DD expansion.
\end{enumerate}

\bibliographystyle{fedcsisIEEEtran}
\bibliography{bibtex/bib/IEEEabrv,dd_sources,jsp_sources}

\end{document}